\documentclass[12pt]{amsart}
\usepackage{graphicx}
\usepackage{epstopdf}
\RequirePackage{color}
\RequirePackage[colorlinks, urlcolor=my-blue,linkcolor=my-red,citecolor=my-green]{hyperref}
\definecolor{my-blue}{rgb}{0.0,0.0,0.6}
\definecolor{my-red}{rgb}{0.5,0.0,0.0}
\definecolor{my-green}{rgb}{0.0,0.5,0.0}
\definecolor{nicos-red}{rgb}{0.75,0.0,0.0}
\definecolor{light-gray}{gray}{0.6}
\definecolor{really-light-gray}{gray}{0.8}
\definecolor{sussexg}{rgb}{0.0,0.5,0.5}
\definecolor{sussexp}{rgb}{0.5,0.0,0.5}

\usepackage{geometry}                
\usepackage[parfill]{parskip}    

\expandafter\let\csname equation*\endcsname\relax
\expandafter\let\csname endequation*\endcsname\relax
\usepackage{amsmath,amssymb,amsthm}

\usepackage{amsfonts}
\usepackage{hyperref}
\usepackage{latexsym}
\usepackage{array}
\usepackage{subfig}

\makeatletter
\newcommand{\addresseshere}{%
  \enddoc@text\let\enddoc@text\relax
}
\makeatother

\begin{document}

\title[Low traffic in a double auction model]{Low-traffic limit and first-passage times for a simple model of the
continuous double auction}

\author{Enrico Scalas}

\address{Department of Mathematics, University of Sussex, Brighton, UK  and  BCAM, Basque Center for Applied Mathematics, Bilbao, Spain }
\email{e.scalas@sussex.ac.uk}

\author{Fabio Rapallo}

\address{DISIT, Universit\`a del Piemonte Orientale, Alessandria, Italy}
\email{fabio.rapallo@uniupo.it}

\author{Tijana Radivojevi\'{c}}

\address{BCAM, Basque Center for Applied Mathematics, Bilbao, Spain}
\email{tradivojevic@bcamath.org}





\begin{abstract}
We consider a simplified model of the continuous double auction where prices are integers varying from $1$ to $N$ with limit orders and market orders, but quantity per order limited to a single share. For this model, the order process is equivalent to two $M/M/1$ queues. We study the behaviour of the auction in the low-traffic limit where limit orders are immediately transformed into market orders. In this limit, the distribution of prices can be computed exactly and gives a reasonable approximation of the price distribution when the ratio between the rate of order arrivals and the rate of order executions is below $1/2$. This is further confirmed by the analysis of the first passage time in $1$ or $N$.
\end{abstract}



\maketitle


\section{Introduction}

Most of the regulated markets in the world implement a trading mechanism known as the {\em continuous double auction} to match supply and demand. This mechanism has two sides. On the supply side there are orders to sell and on the demand side there are orders to buy. Hence, the auction is called {\em double}. Moreover it occurs in continuous time. Hence, it is called {\em continuous}.

In recent years, the theory of this auction has gained more and more interest. In particular, it has been shown that appropriate models of the double auction can be mapped in a multi-class queue \cite{blanchet13}, so that its ergodic properties and the limiting invariant distribution can be studied using established techniques \cite{ferrari09}.

In this paper, we consider a simplified model (see \cite{radivojevic2014} and references therein) where prices take $N$ integer values from $1$ to $N$. Only two types of orders are considered: {\em limit orders} and {\em market orders}. In their turn, limit orders can be either {\color{black}orders to sell} a single share at a price not lower than a given amount ({\em asks}) or {\color{black}orders to buy} a single share at a price not higher than a given amount ({\em bids}). In other words, the quantity attached to every limit order is always $1$. Among all the asks, the {\em best ask} is the smallest ask price, whereas the {\em best bid} is the largest bid price. The best bid is always strictly smaller than the best ask. Market orders have also two sides: either they accept the available best bid or the available best ask.
For the sake of simplicity, limit ask orders and limit bid orders arrive according to a Poisson process at a rate $\lambda_a$ and $\lambda_b$, respectively. In the following, we assume symmetry, i.e. $\lambda_a = \lambda_b = \lambda$. Market orders to buy and market orders to sell arrive separated by durations following the exponential distribution with parameter $\mu_b$ and $\mu_a$, respectively. Again, symmetry is assumed, namely $\mu_a = \mu_b = \mu$.
Limit ask orders follow the uniform distribution in the interval from ${p}_b + 1$ to ${p}_b + n$, where
${p}_b$ is the current best bid and $n \geq 1$ is a parameter of the model. Similarly, limit bid orders are uniformly drawn from the interval
${p}_a - n$ to ${p}_a -1$, where ${p}_a$ is the the current best ask. The accessible states of the auction are limited by the condition ${p}_b < {p}_a$. When ${p}_a$ is between $1$ and $n$ (${p}_b$ between $N-n+1$ and $N$),
the bid (respectively ask) interval is restricted correspondingly. For instance, if ${p}_a = 1$, bids are impossible. The parameter $n$ acts as a cut-off for price jumps.
Eventually, if no orders are present in the auction, the next bid, $b$,
is uniformly chosen in the interval $p-n \leq b \leq p$ and the next ask, $a$, is uniformly taken from $p \leq a \leq p+n$, where $p$ is the price of the last trade.
Specifying an initial price (the opening auction price) is sufficient to start the auction. A short remark is necessary at this stage: It turns out that order inter-arrival times are not exponentially distributed in real markets (see \cite{scalas06b} and references therein). This means, that the above description in terms of $M/M/1$ processes should be replaced by a semi-Markov description in terms of $G/G/1$ processes. However, in this paper, for the sake of simplicity, we will limit our analysis to the Markovian case.

The model described above is essentially the same as in \cite{smith02} and in \cite{zuo09}. It is a
{\em zero intelligence agent-based model} \cite{gode}. As already mentioned in \cite{radivojevic2014}, this version of the model does not use the uniform distribution over $[0, \infty)$ as in \cite{smith02} and it is not limited to the case in which limit orders arrive only at the best bid/ask price as in \cite{zuo09}. A preliminary discussion of this version was presented in \cite{RASold}. This model was extensively studied in \cite{delarrard}, in the case in which price movements equal one tick. These authors also studied the heavy-traffic limit \cite{contdelarrard} where functional limit theorems are available leading to diffusion approximations \cite{kingman, whitt}.

In \cite{radivojevic2014}, the focus was on the ergodic properties of the model. Based on the fact that the order process is equivalent to two independent $M/M/1$ queues, it was shown that
there are three regimes depending on the value of $\rho = \lambda / \mu$. For $0 < \rho < 1$,  prices are free to fluctuate over the full price range and statistical equilibrium is reached (ergodic regime). For $\rho \geq 1$, the auction is in a non-ergodic regime which stabilizes prices. Due to the presence of the parameter $n$, there is an additional transition. If $1 \leq \rho < n$,  prices can still fluctuate in a limited range, whereas for $\rho \geq n$, prices eventually fluctuate between two values. This regime cannot be found if one only considers the case $n=1$.

In the following, we further characterize the ergodic regime by considering the so-called {\em low-traffic} limit where $\rho \ll 1$. It is a limit where analytic results are available for the price dynamics as discussed below. Moreover, we study the first-passage time of the auction in $1$ or in $N$. It turns out that this analysis provides useful approximations for the behaviour of the auction when $\rho < 1/2$.

\section{The low-traffic limit}

In the low-traffic limit ($\rho \ll 1$), when limit orders arrive, they are immediately transformed into market
orders. The book is almost always empty. In this limit, it is possible to explicitly write the transition
probabilities for the price process and study the price Markov chain for any value of $n$ and $N$. To
give an idea on how to proceed, let us assume that the initial price is $p$. Then, the conditional probability of a bid is given by
\begin{equation} \label{form1}
{\mathbb P}(B_1 = b|P_0 = p) = \left\{\begin{array}{lll}
0 & \ \mathrm{ if } \ & b < p - n \ \mathrm{ or } \ b > p \\
\\
\frac 1 p & \ \mathrm{ if } \ & 1 \leq b \leq p \leq n \\
\\
\frac 1 {n + 1} & \ \mathrm{ if } \ & 1 \leq p - n \leq b \leq p.
\end{array} \right.
\end{equation}
This bid is immediately accepted and it becomes the next price. A similar set of equations can be written
for the asks conditioned to the initial price.
\begin{equation}  \label{form2}
{\mathbb P}(A_1 = a|P_0 = p) = \left\{\begin{array}{lll}
0 & \ \mathrm{ if } \  & a < p \ \mathrm{ or } \ a > p + n \\
\\
\frac 1 {N - p + 1} & \ \mathrm{ if } \ & N - n + 1 \leq p \leq a \leq  N \\
\\
\frac 1 {n+1} & \ \mathrm{ if } \ & p \leq a \leq p + n \leq N.
\end{array} \right.
\end{equation}
Both equations are an immediate consequence of the model definition.
For a full characterization of the price Markov chain, the distribution of the initial price is needed. For instance, if the initial price is chosen uniformly, the probability of an initial price is $1/N$; if the chain starts from a given price, the probability of this price is $1$ and the probabilities of all the other prices are $0$, and so on. In a symmetric auction, for which $\lambda_a = \lambda_b = \lambda$ and $\mu_a = \mu_b = \mu$, the probability of a bid arriving is $1/2$ and it is equal to the probability of arrival of an ask. Therefore, in the low traffic regime, the transition probability
for prices is given by
\begin{equation}
\label{transprob}
P_{p,p'}={\mathbb P}(P_1 = p'|P_0 = p) = \frac 1 2 {\mathbb P}(A_1 = p' |P_0 = p) + \frac 1 2 {\mathbb P}(B_1 = p' |P_0 = p).
\end{equation}

Let us now assume that the limit $\rho \ll 1$ is realized by keeping the arrival rate $\lambda$ finite and letting $\mu \gg \lambda$. Then, if $N(t)$ denotes the number
of transactions up to time $t$, we have that $N(t)$ is Poisson distributed with parameter $2\lambda$, given that the auction has two sides (either a bid arrives with rate $\lambda$ or an ask arrives with rate $\lambda$). In fact, $N(t)$ is the superposition of two Poisson processes with
parameter $\lambda$. In other words, the price process can be seen as an embedded Markov chain characterized by the transition probability
\eqref{transprob} subordinated to the Poisson process $N(t)$. Once this remark is made, it is safe to focus on the embedded chain and study its properties. In particular we are interested in the convergence of the price probability. First of all, we notice that after any transactions, the double auction is exactly in the same situation as in the initial case, except for the fact that the price probability varies with time. In other words, the Markov chain defined above is homogeneous. From the study of the transition probability, one can further infer that the Markov chain is irreducible. In fact, it is possible to reach any price from any other price. Moreover, given that the diagonal terms {\color{black} of the Markov transition matrix} are all positive, meaning that there is a finite probability for the price not to change at every step, we can conclude that our Markov chain is aperiodic. Being irreducible and aperiodic, our chain has a unique invariant distribution and this is an equilibrium distribution.

In order to illustrate the above findings, let us consider a specific example with $N=10$ prices and $n=2$. In this case, the price transition probability matrix is
\begin{equation}
\label{matrix}
P = \begin{pmatrix}
	4/6 & 1/6 & 1/6 & 0 & 0& 0& 0& 0& 0 &0 \\
	1/4 & 5/12 & 1/6 & 1/6 & 0& 0& 0& 0& 0 &0 \\
	1/6 & 1/6 & 1/3 & 1/6 & 1/6& 0& 0& 0& 0 &0 \\
	0 & 1/6 & 1/6 & 1/3 & 1/6& 1/6& 0& 0& 0 &0 \\
	0 & 0 & 1/6 & 1/6 & 1/3& 1/6& 1/6& 0& 0 &0 \\
	0 & 0 & 0 & 1/6 & 1/6& 1/3& 1/6& 1/6& 0 &0 \\
	0 & 0 & 0 & 0 & 1/6& 1/6& 1/3& 1/6& 1/6 &0 \\
	0 & 0 & 0 & 0 & 0& 1/6& 1/6& 1/3& 1/6 &1/6 \\
	0 & 0 & 0 & 0 & 0& 0& 1/6& 1/6& 5/12 &1/4 \\
	0 & 0 & 0 & 0 & 0& 0& 0& 1/6& 1/6 &4/6 \\
\end{pmatrix}.
\end{equation}
The invariant distribution is obtained by looking for the left eigenvector with unit eigenvalue, namely
\begin{equation}
\label{eig}
 \pi P= \pi,
\end{equation}
which, in the case of \eqref{matrix} gives
\begin{equation}
\pi = \left(0.1171, 0.0895,  0.1,  0.0961,  0.0974, 0.0974,  0.0961,  0.1,  0.0895,  0.1171 \right).
\end{equation}
In the Appendix A, we present a general algorithm to find the invariant distribution of prices. For $n=1$, there is a remarkable result. In fact, in this case, the transition matrix is a symmetric, doubly-stochastic matrix. Since $P\mathbf1 = \mathbf1$ and $\mathbf1^TP = \mathbf1^T$ (because row sums and column sums are 1), then
\begin{equation}
\label{uniform}
\frac{1}{N}\mathbf1^TP = \frac{1}{N}\mathbf1^T,
\end{equation}
and the uniform distribution is the invariant distribution for the Markov chain.

\begin{figure}
\includegraphics[width=90mm,height=80mm]{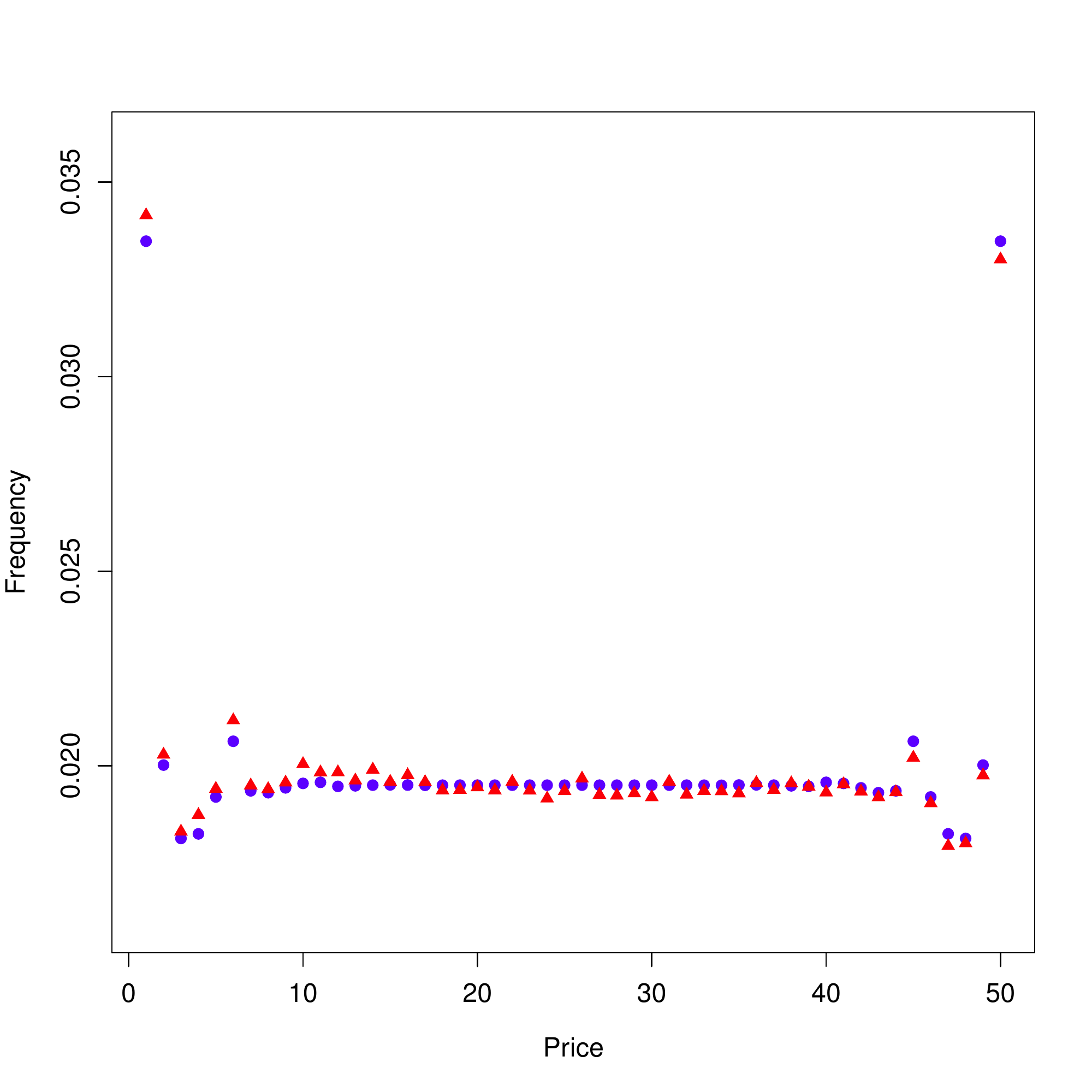}
\caption{Equilibrium price distribution in the low-traffic limit in the case $N=50$, $n=5$ (circles). The triangles denote the price frequency for a Monte Carlo simulation of the double auction with $\rho=10^{-4}$ after $10^6$ steps.}
\label{fig1}
\end{figure}

\begin{figure}
\includegraphics[width=90mm,height=80mm]{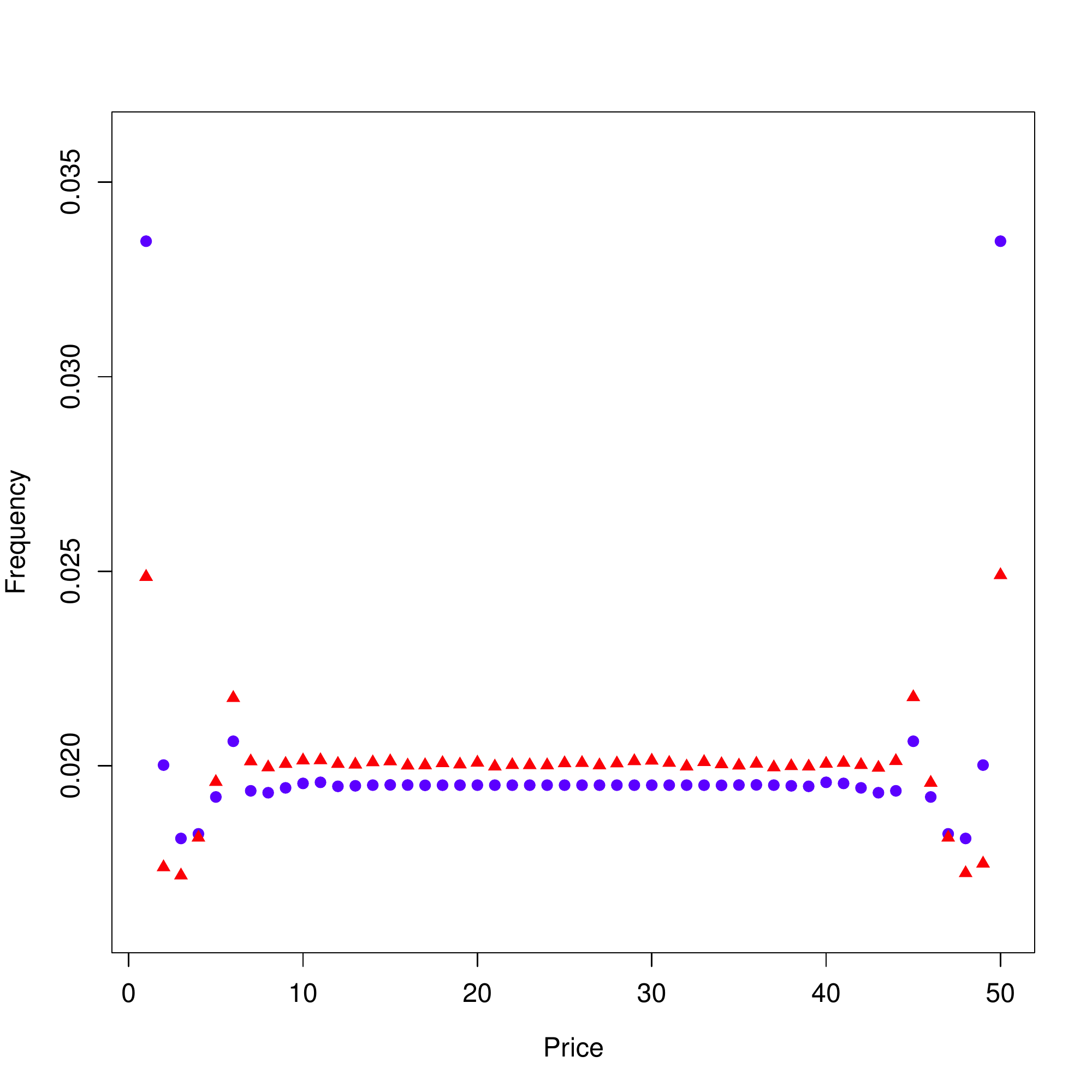}
\caption{Comparison between the equilibrium price distribution in the low-traffic limit in the case $N=50$, $n=5$ (circles) and a a Monte Carlo simulation (triangles) of the double auction with $\rho=0.3$ after $10^9$ steps.}
\label{fig2}
\end{figure}

\begin{figure}
\includegraphics[width=90mm,height=80mm]{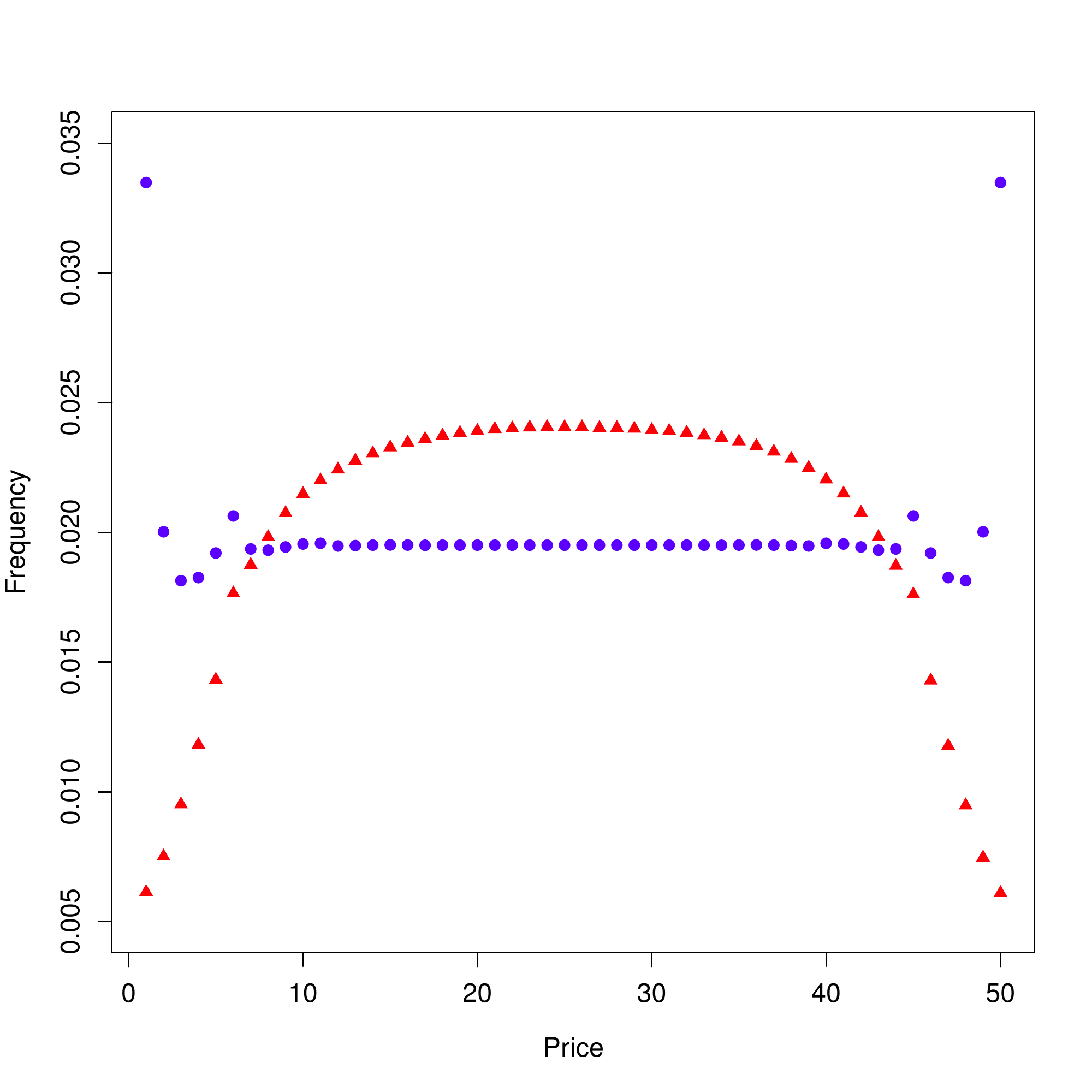}
\caption{Comparison between the equilibrium price distribution in the low-traffic limit in the case $N=50$, $n=5$ (circles) and a a Monte Carlo simulation (triangles) of the double auction with $\rho=0.9$ after $10^9$ steps.}
\label{fig3}
\end{figure}

Always for purpose of illustration, in Figure \ref{fig1}, we plot the low-traffic limit price distribution for the case $N=50$, $n=5$ and we compare it with the frequency with which states appear after equilibration in a Monte Carlo simulation of the chain after $10^6$ iterations for $\rho=10^{-4}$. In this case $10^6$ iterations are already sufficient to show that the agreement between the low-traffic-limit approximation and the result of Monte Carlo simulations is good.
It is striking to observe that the price distribution in the low-traffic limit is still a reasonable approximation when $\rho=0.3$ as shown in Figure \ref{fig2}.
The approximation breaks down for $\rho \geq 0.5$ as shown in Figure \ref{fig3} for the case $\rho=0.9$.
This behaviour of the price distributions leads to a different behaviour for first passage times at the boundary prices. In fact, for $\rho > 0.5$ and $n>1$ the residence time of the systems close to the boundaries becomes negligible as shown by Figure \ref{fig3}, leading to an increase of the value of the average first passage time.
%
%
%
%

\section{First passage times}

In this section we shall focus on first passage times. Given a double auction with $N$ possible prices, labeled with the integers $1, \ldots, N$, we fix the initial price at the median point $\lfloor \frac {N+1} 2 \rfloor$, where $\lfloor \cdot \rfloor$ denotes the floor operator. We study the random variable $T$: The first passage time at $1$ or at $N$; for this reason, our problem belongs to the class of two-barrier problems, and, for simplicity, we will assume that the price walk is symmetric. The behaviour of the first passage time distribution has been studied in several similar problems, both with theoretical results on the exact or asymptotic distribution, and through simulation studies, see e.g. \cite{bondesson} and \cite{antalredner}.

Here, we present the results of a simulation study to investigate the main features of the distribution of $\log(T)$. In particular, we compare such a distribution with the theoretical distribution derived under the low traffic assumption. Since the distribution of $T$ is highly skewed with a fat upper tail, as shown in figures, all the plots reported here refer to the distribution of its natural logarithm $\log(T)$. For sake of simplicity, we have performed all comparisons in this section with the parameter $\mu = 1$ fixed. With this assumption, we get the low traffic limit if $\rho = \lambda \ll 1$.

\begin{figure}
\includegraphics[width=90mm]{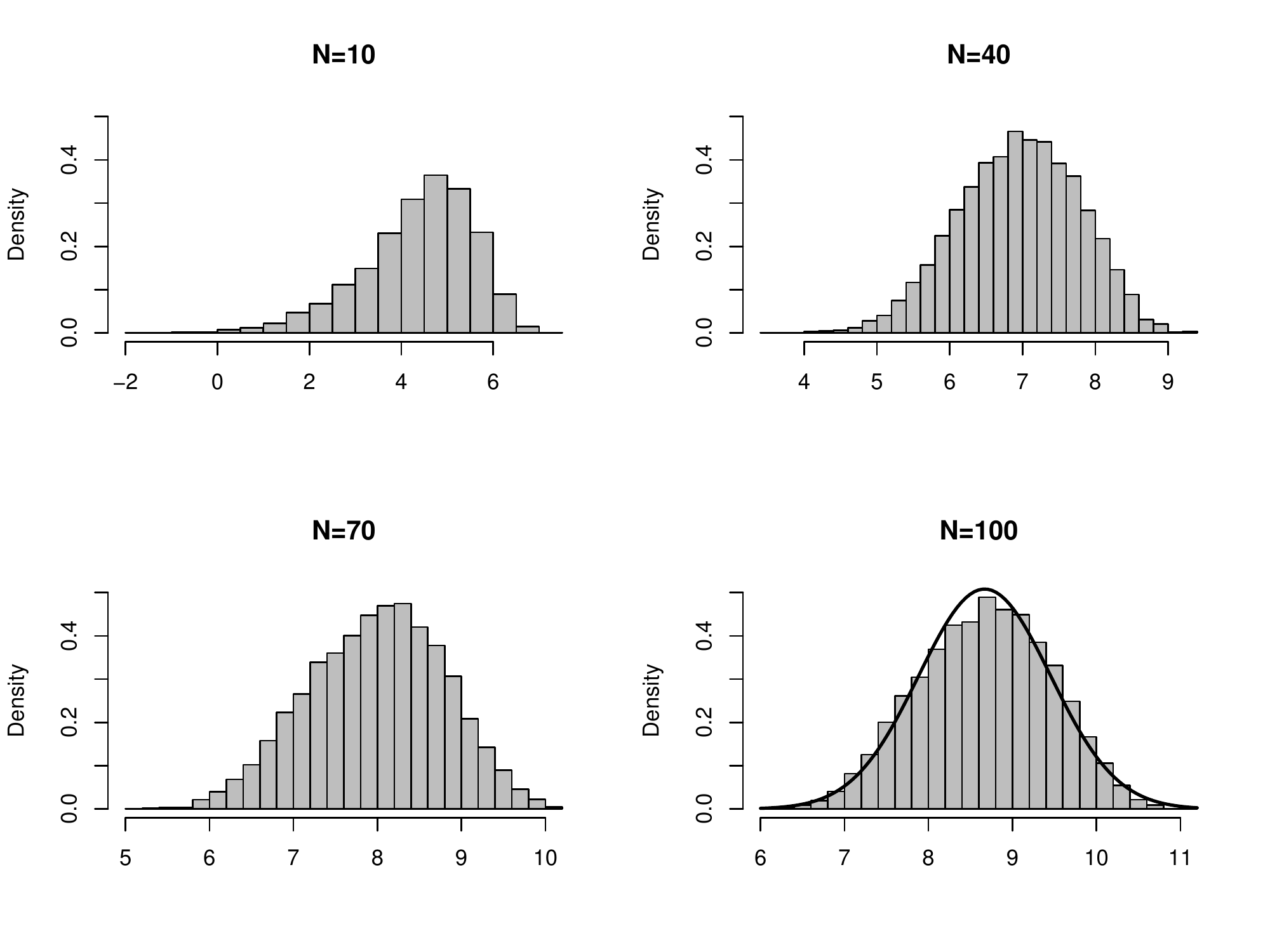}
\caption{Distribution of $\log(T)$ for $n=5$ fixed and $\rho=0.02$.}
\label{fig-hist1}
\end{figure}

In Figure \ref{fig-hist1} the histograms of $\log(T)$ for $n=5$ and for $\rho=0.02$ are displayed for $4$ different values of $N$, namely $N=10,40,70,100$. One can observe that the shape of the distribution is skewed for small values of $N$, while it approaches a Gaussian distribution in the case $N=100$ (the best-fit normal curve is plotted together with the histogram). The simulations have been implemented in R \cite{erre}, and all histograms in this study are based on $10,000$ Monte-Carlo replicates.

\begin{figure}
\includegraphics[width=90mm]{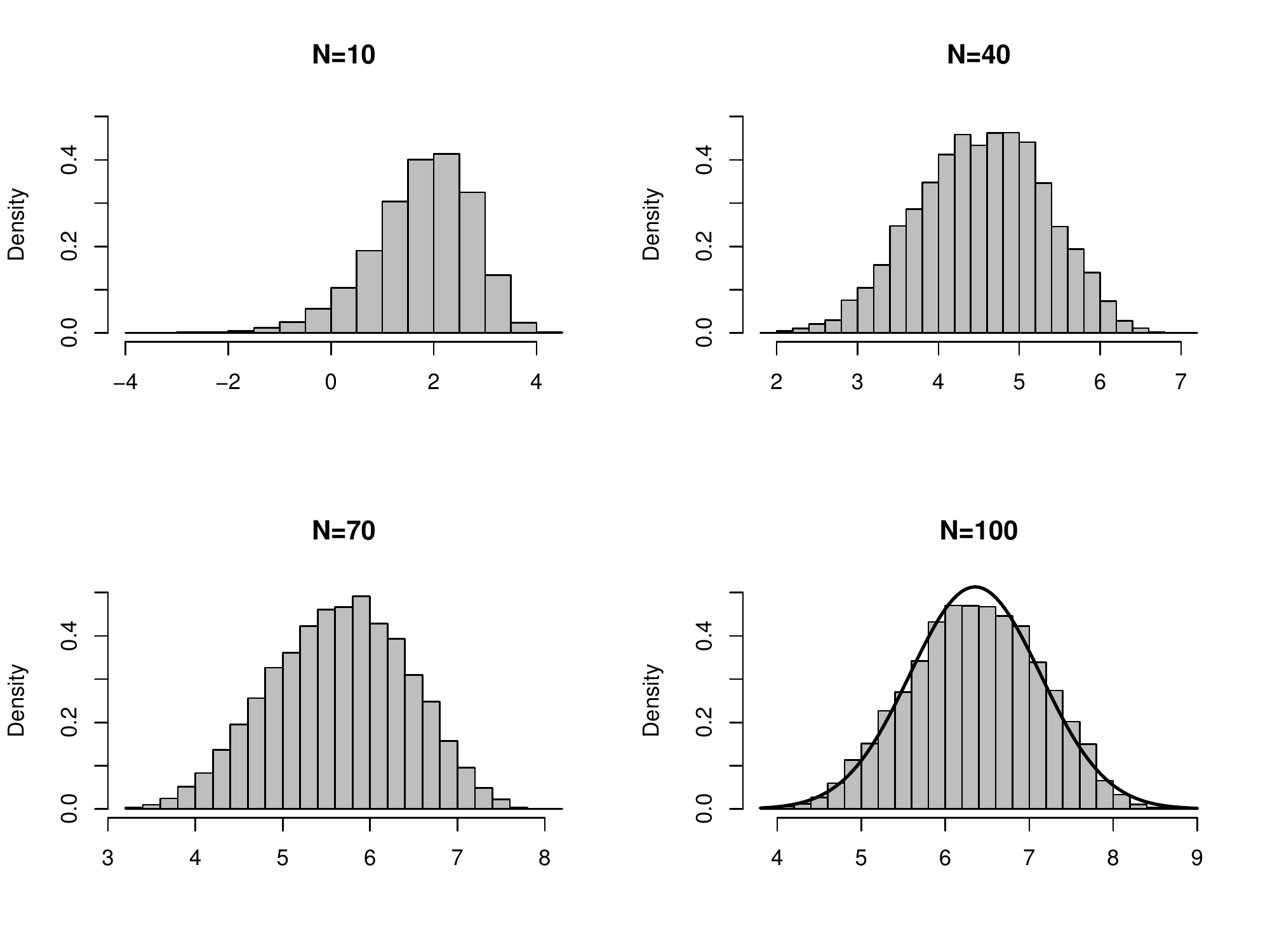}
\caption{Distribution of $\log(T)$ for for $n=5$ fixed and $\rho=0.5$.}
\label{fig-hist2}
\end{figure}

Figure \ref{fig-hist2} refers to the same settings as above, but with the ratio fixed at $\rho=0.5$. The distribution of $\log(T)$ has almost the same shape for $\rho=0.5$ and for $\rho=0.02$. In agreement with the conclusions of the previous section, this suggests that the behaviour of the low traffic limit is a good approximation also for values of $\rho$ up to $0.5$ also in terms of the first passage time distribution.

\begin{figure}
\includegraphics[width=90mm]{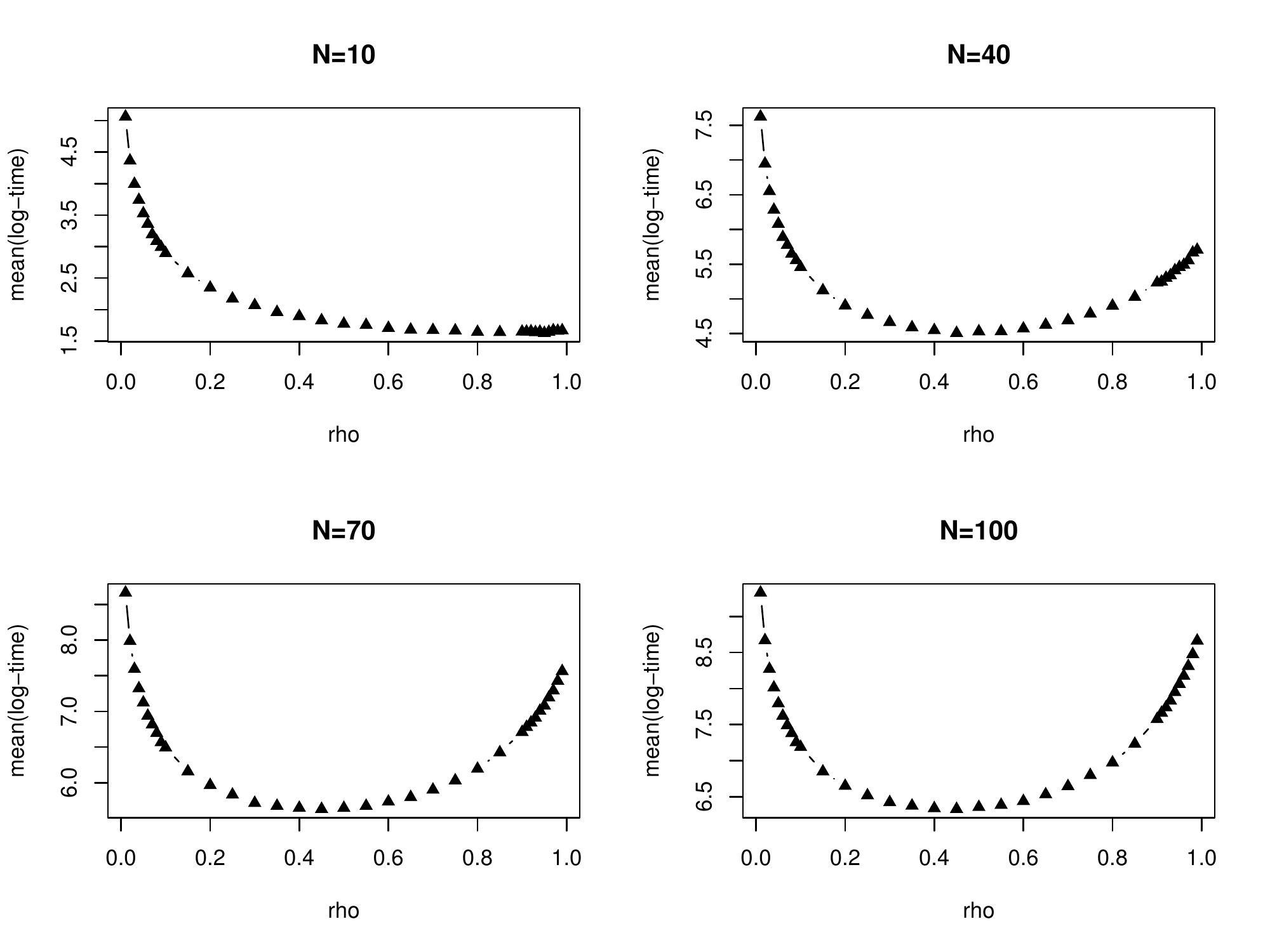}
\caption{Mean of $\log(T)$ for $n=5$ and $N=10,40,70,100$.}
\label{fig-means}
\end{figure}

Finally, in Figure \ref{fig-means}, the means of $\log(T)$ as a function of $\rho$ under various choices of $N$ ($n=5$ fixed) are displayed. A minimum occur between $\rho=0.4$ and $\rho=0.5$ except for the first experimental setting ($N=10$).

To complete the simulation study, we have compared the distribution of $T$ with an approximation suggested by the results of the previous section. Basically, we adapt here a known formula for a discrete two-barrier problem. Such formula gives the distribution of the number of price changes needed to reach the boundary. Then, the parameter $\rho$ controls the proportion of orders leading to a price change, and therefore it defines the discrete time distribution of the number of orders needed to reach the boundary. Remember that in the low traffic limit when limit orders arrive, they are immediately transformed into market orders. Finally, we add suitable exponential distributions to switch to the continuous framework of our model. To avoid problems in some formulae, we assume here $N$ to be odd, so that $(N+1)/2$ is always integer. We have limited our study to the case $n=1$ in order to avoid further technicalities in the formulae and to capture the major features of the model. The approximation is built up in three steps, as detailed below:
\begin{itemize}
\item First, consider a discrete first passage time ${\widetilde T}^{(d)}$ in a simple symmetric random walk with two reflecting barriers and discrete $\pm 1$ steps, whose distribution is
\begin{equation}
\label{discr-formula}
{\mathbb P}({\widetilde T}^{(d)} = {\tilde h}) = \frac {2} {N-1} \sum_{k=1}^{N-2} (-1)^{k+1} \sin \left( \frac {k \pi}{N-1} \right) \cos^{{\tilde h}-1 } \left( \frac{k \pi} {N-1} \right) \sin \left( \frac {k \pi} 2 \right)
\end{equation}
for $\tilde h \geq 1$. The above distribution can be found in \cite{fellerI} and is extensively discussed with several generalizations in \cite{carlsund}.

\item The rate of arrival of limit orders over all orders is $\lambda / (\lambda + \mu) = \rho/(\rho+1)$, and in the low traffic approximation all limit orders arrive when the book is empty. A limit order to buy (resp. to sell) fixes the price at the old price $p$ or at $p+1$ (resp. $p-1$) with probability $1/2$ each. Therefore, at any given time, the price changes with rate $\frac {\rho} {2(\rho + 1)}$. Thus, given ${\widetilde T}^{(d)}={\tilde h}$, consider a Negative Binomial variable $NB$ with parameters $\tilde h$ and $\frac {\rho} {2(\rho + 1)}$, and define $T^{(d)}=NB+1$. The variable $T^{(d)}$ is again a discrete random variable and it counts the number of events after the change of price is actually performed;

\item The interarrival time between two consecutive events follows an exponential distribution with mean $2\mu(1+\rho)$. Therefore, given $T^{(d)} = h$, the first passage time is approximated by a random variable following a Gamma distribution with parameters $h$ and $\frac 1 {2\mu(1+\rho)}$.
\end{itemize}
In conclusion, the distribution of the first passage time can be approximated by a suitable mixture $T^{(a)}$ of Gamma distributions, whose parameters are computed according to the formula in Eq. \eqref{discr-formula} for the discrete case. Notice that in the previous construction the low-traffic hypothesis is used only when we assume the book empty when a limit order arrives.

\begin{figure}
\includegraphics[width=90mm]{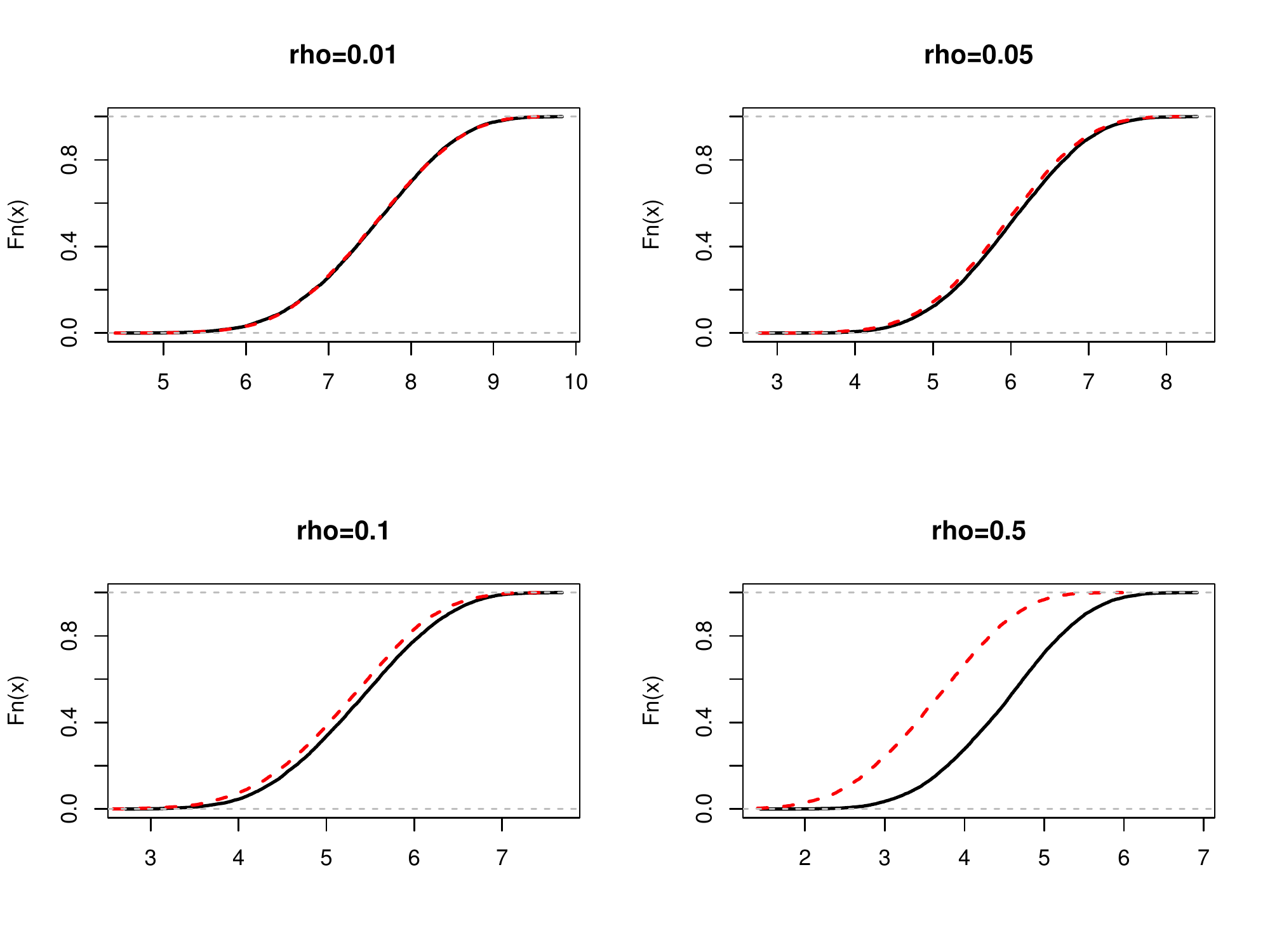}
\caption{ECDFs of $\log(T)$ (in black) and its low traffic approximation $\log(T^{(a)})$ (dashed, in red) for $N=11$, $n=1$.}
\label{fig-ecdfs}
\end{figure}

To show that this approximation works well for small values of $\rho$, we have plotted in Figure \ref{fig-ecdfs} the empirical cumulative distribution function (ECDF) of $\log(T)$ and (a Monte-Carlo approximation of) the distribution of $\log(T^{(a)})$ for $4$ different values of $\rho$ in the case $N=11$. Also in these simulations we have considered only the case $\mu = 1$.

We can observe in Figure \ref{fig-ecdfs} that for $\rho=0.01$, $\rho=0.05$ the simulated distribution and its theoretical approximation are nearly identical (the $p$-value of the Kolmogorov-Smirnov based on 10,000 Monte Carlo replicates is $0.7212$ for $\rho=0.01$ and $0.0116$ for $\rho=0.05$. When $\rho=0.1$ the two distributions show some discrepancies, while in the last scenario ($\rho=0.5$) the approximation fails. The low traffic approximation $T^{(a)}$ tends to underestimate the distribution of $T$. This behaviour is observed also in the $\rho=0.1$ case, but it is clearer in the $\rho=0.5$ case, as expected.

When $n>1$ a formula like that in Eq. \eqref{discr-formula} is no longer available. However, we can analyze the low-traffic approximation by studying the expected values of the distributions. In fact, the expected value of the first passage time $\mu_T^{(d)}$ of the discrete chain in the low-traffic approximation can be computed through the linear system
\begin{equation}
(I-P_{2,N-1})x = 1
\end{equation}
where $P_{2,N-1}$ is the transition matrix restricted to the transient states, $I$ is the $(N-2)\times (N-2)$ identity matrix, and $1$ is a column vector of 1 with dimension $N-2$ (see e.g. \cite{kemenysnell} for details). Then, the mean time in the continuous setting $\mu_T$ is obtained by scaling $\mu_T^{(d)}$ by a factor $1/(2\rho)$, following the same reasoning as above. In Table \ref{tabmeantimes}, the means $\overline T$ of the Monte Carlo simulations and the theoretical expected value under the low-traffic approximation $\mu_T$ are given for several settings. For $\rho$ up to $0.05$ the approximation works well, and the relative error is less than $10\%$ in all settings, while for $\rho=0.1$ and $\rho=0.5$ the differences become relevant, especially in the latter case.

\begin{table}
\begin{center}
\begin{tabular}{cc|ccc|ccc|ccc}
 & & & $\rho=0.01$ & &  & $\rho=0.02$ & &  & $\rho=0.05$ &  \\
$N$ & $n$ & $\overline T$ & $\mu_T$ & $\Delta\%$ & $\overline T$ & $\mu_T$ & $\Delta\%$ & $\overline T$ & $\mu_T$ & $\Delta\%$ \\ \hline
10 & 5 & 273.17 & 275.67 & $-0.91$ & 136.58 & 137.76 & $-0.85$ &  54.63 & 56.44 & $-3.20$ \\
40 & 5 & 2787.28 & 2821.84 & $-1.22$ & 1393.64 & 1433.86 & $ -2.80$ &  557.46 & 598.03 & $-6.78$ \\
40 & 10 & 1176.20 & 1197.19 & $-1.75$ & 588.10 & 604.54 & $-2.72$ &  235.24 & 248.15 & $-5.20$ \\
80 & 5 & 9939.48 & 9888.86 & $+0.51$ & 4969.74 & 5181.18 & $-4.08$ &  1987.90 & 2164.12 & $ -8.14$ \\
80 & 20 & 1598.81 & 1617.58 & $-1.16$ & 799.41 & 824.12 & $-3.00$ &  319.76 & 342.19 & $ -6.55$ \\
100 & 5 & 15151.93 & 15183.58 & $-0.21$ &  7575.97 & 7807.49 & $-2.97$ &  3030.39 & 3267.66 & $ -7.26$ \\
100 & 25 & 1763.36 & 1768.67 & $ -0.30$ &   881.68 & 903.18 & $-2.38$ & 352.67 & 369.58 & $  -4.58$ \\
\hline
\end{tabular}
\medskip
\begin{tabular}{cc|ccc|ccc}
 &  & & $\rho=0.10$ & &  & $\rho=0.50$ &  \\
$N$ & $n$ & $\overline T$ & $\mu_T$ & $\Delta\%$ & $\overline T$ & $\mu_T$ & $\Delta\%$ \\ \hline
10 & 5 & 27.32 & 29.47 & $-7.31$ & 5.46 & 8.72 & $-37.35$ \\
40 & 5 & 278.73 & 319.18 & $-12.67$ & 55.75 &124.99 & $-55.40$ \\
40 & 10 & 117.62 & 132.71 & $-11.37$ & 23.52 & 45.54 & $-48.35$ \\
80 & 5 & 993.95 & 1159.32 & $-14.26$ &  198.79 & 495.32 & $ -59.87$ \\
80 & 20 & 159.88 & 182.86 & $-12.57$ & 31.98 &  62.29 & $  -48.66$ \\
100 & 5  & 1515.19  & 1780.02  & $-14.88$ & 303.04 &  769.96 & $  -60.64$ \\
100 & 25 & 176.34  & 198.62   & $-11.22$ & 35.27 &  70.41 & $  -49.91$ \\
\hline
\end{tabular}
\end{center}
\caption{Average times in the real settings and the corresponding low-traffic approximation for different values of $N$, $n$ and $\rho$.}
\label{tabmeantimes}
\end{table}

\section{Summary and conclusions}

In this paper, we characterized the ergodic regime of a simple model for the continuous double auction in the low-traffic limit $\rho \ll 1$. In this limit, the price distribution can be derived for any value of the model parameters $n$ and $N$. Explicit numerical procedures to find the price distributions are given in the Supplemental Material. We also showed that these results give a reasonable approximation of the auction behaviour for $\rho < 1/2$. We further studied the first passage time $T$ in $1$ or $N$ using Monte Carlo simulations. We noticed that the low-traffic limit approximation for $T$ works reasonably for $\rho \ll 1/2$ in this case.

There are several open questions we would like to answer. A natural extension of this simple model is its semi-Markov version in which non exponential distributions for waiting times between events are introduced. In such an extension, the behavior of the embedded chain does not change, but the mixing time of the chain changes. A particularly interesting case is when the distribution of waiting times is heavy-tailed with infinite mean. This is linked to recent results of ours on semi-Markov graph dynamics \cite{raberto11, georgiou15}.
A further research direction worth exploring is considering non-independent processes for limit and market orders.

\section*{Acknowledgements}

T. Radivojevi\'c and E. Scalas wish to thank J. Anselmi for useful discussion. E. Scalas thanks N. Georgiou for discussion on multi-class queues.


\addresseshere

\newpage

\appendix
\section*{Appendix A}

In this appendix, we present a general algorithm to find the invariant distribution for the price Markov chain. First, we observe that, in general, the transition probability matrix is a stochastic block tri-diagonal (symmetric in the inner part) matrix of the form:
\begin{equation}
P = \begin{pmatrix}
	\color{green}{D_0} & \color{red}{A} &  0 & \dots & &0 \\
	\color{red}{A^T} & \color{blue}{D} & \color{red}A & \dots & &0 \\
	0 & \color{red}A^T & \color{blue}D &  \dots & &0 \\
	\vdots &  & &  &\color{blue} D &\color{red}A \\
	0 &  & \dots & & \color{red}A^T &\color{green}\widetilde{D_0} \\
\end{pmatrix},
\end{equation}
where $A^T$ is the transpose of the block $A$ and $\widetilde{D_0}=(d_{n-i+1,n-j+1})$ for $d_{ij}$ being elements of the block $D_0$.
There is some freedom in the choice of the blocks $D_0$, $D$ and $A$ and their transformations whose dimensions depend on the value of $n$. In the case of the transition probability (4) in the paper, this is a possible choice of blocks:
\begin{equation}
P = \begin{pmatrix}
	\color{green}{4/6} & \color{green}{1/6} & \color{red}{1/6} & \color{red}{0} & 0& 0& 0& 0& 0 &0 \\
	\color{green}{1/4} & \color{green}{5/12} & \color{red}{1/6} & \color{red}{1/6} & 0& 0& 0& 0& 0 &0 \\
	\color{red}{1/6} & \color{red}{1/6} & \color{blue}{1/3} & \color{blue}{1/6} & \color{red}1/6& \color{red}0& 0& 0& 0 &0 \\
	\color{red}{0} & \color{red}{1/6} & \color{blue}{1/6} & \color{blue}{1/3} & \color{red}1/6& \color{red}1/6& 0& 0& 0 &0 \\
	0 & 0 & \color{red}1/6 & \color{red}1/6 & \color{blue}1/3& \color{blue}1/6& \color{red}1/6& \color{red}0& 0 &0 \\
	0 & 0 & \color{red}0 & \color{red}1/6 & \color{blue}1/6& \color{blue}1/3& \color{red}1/6& \color{red}1/6& 0 &0 \\
	0 & 0 & 0 & 0 & \color{red}1/6& \color{red}1/6& \color{blue}1/3& \color{blue}1/6& \color{red}1/6 &\color{red}0 \\
	0 & 0 & 0 & 0 & \color{red}0& \color{red}1/6& \color{blue}1/6& \color{blue}1/3& \color{red}1/6 &\color{red}1/6 \\
	0 & 0 & 0 & 0 & 0& 0& \color{red}1/6& \color{red}1/6& \color{green}5/12 &\color{green}1/4 \\
	0 & 0 & 0 & 0 & 0& 0& \color{red}0& \color{red}1/6& \color{green}1/6 &\color{green}4/6 \\
\end{pmatrix},
\end{equation}
but the following block choice seems to be more convenient
\begin{equation}
P = \begin{pmatrix}
	\color{green}{4/6} & \color{green}{1/6} & \color{red}{1/6} & 0 & 0& 0& 0& 0& 0 &0 \\
	\color{green}{1/4} & \color{green}{5/12} & \color{red}{1/6} & \color{red}1/6 & 0& 0& 0& 0& 0 &0 \\
	\color{red}{1/6} & \color{red}{1/6} & \color{blue}{1/3} & \color{red}1/6 & \color{red}1/6& 0& 0& 0& 0 &0 \\
	0 & \color{red}1/6 & \color{red}1/6 & \color{blue}{1/3} & \color{red}1/6& \color{red}1/6& 0& 0& 0 &0 \\
	0 & 0 & \color{red}1/6 &\color{red} 1/6 & \color{blue}1/3& \color{red}1/6& \color{red}1/6& 0& 0 &0 \\
	0 & 0 & 0 & \color{red}1/6 & \color{red}1/6& \color{blue}1/3& \color{red}1/6& \color{red}1/6& 0 &0 \\
	0 & 0 & 0 & 0 & \color{red}1/6& \color{red}1/6& \color{blue}1/3& \color{red}1/6& \color{red}1/6 &0 \\
	0 & 0 & 0 & 0 & 0& \color{red}1/6& \color{red}1/6& \color{blue}1/3& \color{red}1/6 &\color{red}1/6 \\
	0 & 0 & 0 & 0 & 0& 0& \color{red}1/6& \color{red}1/6& \color{green}5/12 &\color{green}1/4 \\
	0 & 0 & 0 & 0 & 0& 0& 0& \color{red}1/6& \color{red}1/6 &\color{green}4/6 \\
\end{pmatrix}.
\end{equation}
To better see how this generalizes, let us consider the structure of the matrix for $N=10$ and $n=2$, once more
\begin{equation}
P = \begin{pmatrix}
	\color{green}{d_1} & \color{red}{a} & \color{red}{a} & 0 & 0& 0& 0& 0& 0 &0 \\
	\color{green}{b_2} & \color{green}{d_2} & \color{red}{a} & \color{red}{a} & 0& 0& 0& 0& 0 &0 \\
	\color{red}{a} & \color{red}{a} & \color{blue}{2a} & \color{red}{a} & \color{red}{a}& 0& 0& 0& 0 &0 \\
	0 & \color{red}{a} & \color{red}{a} & \color{blue}{2a} & \color{red}{a}& \color{red}{a}& 0& 0& 0 &0 \\
	0 & 0 &  \color{red}{a}& \color{red}{a} & \color{blue}{2a}& \color{red}{a}& \color{red}{a}& 0& 0 &0 \\
	0 & 0 & 0 & \color{red}{a} & \color{red}{a}& \color{blue}{2a}& \color{red}{a}& \color{red}{a}& 0 &0 \\
	0 & 0 & 0 & 0 & \color{red}{a}& \color{red}{a}& \color{blue}{2a}& \color{red}{a}& \color{red}{a} &0 \\
	0 & 0 & 0 & 0 & 0& \color{red}{a}& \color{red}{a}& \color{blue}{2a}& \color{red}{a} &\color{red}{a} \\
	0 & 0 & 0 & 0 & 0& 0& \color{red}{a}& \color{red}{a}& \color{green}{d_2} &\color{green}{b_2} \\
	0 & 0 & 0 & 0 & 0& 0& 0& \color{red}{a}& \color{red}{a} &\color{green}{d_1} \\
\end{pmatrix}.
\end{equation}
For $N=10, n=3$, instead, we have
\begin{equation}
P = \begin{pmatrix}
	\color{green}{d_1} & \color{red}{a} & \color{red}{a} &  \color{red}{a} & 0& 0& 0& 0& 0 &0 \\
	\color{green}{b_2} & \color{green}{d_2} & \color{red}{a} & \color{red}{a} &  \color{red}{a}& 0& 0& 0& 0 &0 \\
	\color{green}{b_3} & \color{green}{b_3} & \color{green}{d_3} & \color{red}{a} & \color{red}{a}& \color{red}{a}& 0& 0& 0 &0 \\
	 \color{red}{a} &  \color{red}{a} &  \color{red}{a} & \color{blue}{2a} &  \color{red}{a}&  \color{red}{a}&  \color{red}{a}& 0& 0 &0 \\
	0 &  \color{red}{a} &  \color{red}{a} &  \color{red}{a} &  \color{blue}{2a}& \color{red}{a}& \color{red}{a}& \color{red}{a}& 0 &0 \\
	0 & 0 & \color{red}{a} & \color{red}{a} & \color{red}{a}& \color{blue}{2a}& \color{red}{a}& \color{red}{a}& \color{red}{a} &0 \\
	0 & 0 & 0 & \color{red}{a} & \color{red}{a}& \color{red}{a}& \color{blue}{2a}& \color{red}{a}& \color{red}{a} &\color{red}{a} \\
	0 & 0 & 0 & 0 & \color{red}{a}& \color{red}{a}& \color{red}{a}& \color{green}{d_3}& \color{green}{b_3} &\color{green}{b_3} \\
	0 & 0 & 0 & 0 & 0& \color{red}{a}& \color{red}{a}& \color{red}{a}& \color{green}{d_2} &\color{green}{b_2} \\
	0 & 0 & 0 & 0 & 0& 0& \color{red}{a}& \color{red}{a}& \color{red}{a} &\color{green}{d_1} \\
\end{pmatrix}
\end{equation}
and so on, where:
 $$a=\frac{1}{2(n+1)} $$
  \begin{equation}
 d_i=\frac{1}{2i}+a
  \end{equation}
$$ b_i=\frac{1}{2i}$$
as a consequence of (3) in the paper and
$$d_1+na=1 $$
 $$b_2+d_2+na=1$$
 $$b_3+b_3+d_3+na=1$$
 \begin{equation}
  \vdots
  \end{equation} 	
$$ (n-1)b_n +d_n +na=1$$
$$ 2na+2a=1$$
as a consequence of the properties of the transition matrix.

The linear system of equations whose solution is the invariant distribution from Eq. (5) in the paper is:
\begin{equation}
\label{system}
\begin{cases}
d_1\pi_1+b_2\pi_2+\dots+b_n\pi_n+a\pi_{n+1}&=\pi_1 \\
a\pi_1+d_2\pi_2+b_3\pi_3+\dots+b_n\pi_n+a\pi_{n+1}+a\pi_{n+2}&=\pi_2 \\
\vdots & \vdots \\
a\pi_{1}+\dots+a\pi_{n-1}+d_n\pi_n+a\pi_{n+1}+\dots+a\pi_{2n}&=\pi_n \\
a\pi_1+\dots+a\pi_n+2a\pi_{n+1}+a\pi_{n+2}\dots+a\pi_{2n+1}&=\pi_{n+1} \\
a\pi_2+\dots+a\pi_{n+1}+2a\pi_{n+2}+a\pi_{n+3}\dots+a\pi_{2n+2}&=\pi_{n+2} \\
\vdots & \vdots \\
a\pi_{N-2n}+\dots+a\pi_{N-n-1}+2a\pi_{N-n}+a\pi_{N-(n-1)}\dots+a\pi_{N}&=\pi_{N-n} \\
a\pi_{N-2n-1}+\dots+a\pi_{N-n}+d_n\pi_{N-(n-1)}+a\pi_{N-(n-2)}+\dots+a\pi_{N}&=\pi_{N-(n-1)} \\
\vdots  & \vdots \\
a\pi_{N-n}+b_{n}\pi_{N-(n-1)}\dots+b_2\pi_{N-1}+d_{1}\pi_{N}&=\pi_{N} \\
\end{cases}
\end{equation}
with the additional equation:
\begin{equation}
\sum_{i=1}^N \pi_i=1.
\end{equation}

A possible numerical \textbf{solver} for the system \eqref{system} given $N$ and $n$, and written in MATLAB is:
{\small
\begin{verbatim}
function prob = pricedistrlt(N,n)
a = 1/2/(n+1);
% defining the transition matrix
d = 2*a*ones(1,N); % blue diagonal values
d(1:n) = 0.5./[1:n] + a; % green  diagonal values
P = diag(d);
for i=1:n % red
    P = P + diag(a*ones(1,N-i),i) + diag(a*ones(1,N-i),-i);
end
for i=2:n % green
    for j=1:i-1
        P(i,j) = 1/2/i;
    end
end
P(N-n+1:N,N-n+1:N) = rot90(P(1:n,1:n),2);

X = sym('x',[1 N]);
x = solve(X*P - X,sum(X)-1);
x = struct2cell(x);
prob = zeros(1,N);
for i=1:N
    prob(i) = x{i};
end
\end{verbatim}
}

The above function requires the Symbolic Math Toolbox in MATLAB. Another possibility is to solve  \eqref{system} using the Matlab function \texttt{eig} after defining the transition matrix in the following way:

{\small
\begin{verbatim}
[prob,l] = eig(P.');
if n==1
    prob = prob(:,end)./sum(prob(:,end));
else
    prob = prob(:,1)./sum(prob(:,1));
end
\end{verbatim}
}

This option is much faster, but less accurate.

\end{document}